\documentclass[12pt,oneside]{amsart}
\usepackage{ifpdf}

\ifpdf
  \usepackage[pdftex]{epsfig}
  \usepackage[pdftex,colorlinks=true,%
              pdftitle={Heuristics on rational points},%
              pdfauthor={Michael Stoll}]{hyperref}
  \newcommand{\Gr}[2]{\psfig{file=#1.pdf,width=#2}}
\else
  \usepackage{epsfig}
  \newcommand{\Gr}[2]{\psfig{file=#1.eps,width=#2}}
\fi

\usepackage{amssymb}  
\usepackage{latexsym} 
\usepackage{comment}

\usepackage[all]{xy}

\DeclareFontEncoding{OT2}{}{} 

\newcommand{\Z}{{\mathbb Z}}
\newcommand{\Q}{{\mathbb Q}}
\newcommand{\R}{{\mathbb R}}

\newcommand{\F}{{\mathbb F}}
\newcommand{\E}{{\mathbb E}}
\newcommand{\BP}{{\mathbb P}}
\newcommand{\BA}{{\mathbb A}}
\newcommand{\To}{\longrightarrow}

\newcommand{\vol}{\operatorname{vol}}

\newcommand{\Prob}{\operatorname{Prob}}

\newcommand{\CC}{{\mathcal C}}
\newcommand{\CD}{{\mathcal D}}
\newcommand{\eps}{\varepsilon}
\newcommand{\sump}{\mathop{\sum\nolimits'}}
\newcommand{\bfx}{\mathbf{x}}
\newcommand{\bfa}{\mathbf{a}}
\newcommand{\bfal}{\boldsymbol{\alpha}}
\newcommand{\bfe}{\mathbf{e}}
\newcommand{\GL}{\operatorname{GL}}
\newcommand{\Gm}{{\mathbb G}_{\text{\rm m}}}

                      {\hspace*{\fill}\nobreak$\Box$\par}

   {\begin{list}{}{\settowidth{\labelwidth}{#1}%
                   \setlength{\leftmargin}{\labelwidth}%
                   \addtolength{\leftmargin}{\labelsep}}}%
   {\end{list}}

\newtheorem{Theorem}{Theorem}
\newtheorem{Lemma}[Theorem]{Lemma}
\newtheorem{Proposition}[Theorem]{Proposition}
\newtheorem{Corollary}[Theorem]{Corollary}
\newtheorem{Conjecture}[Theorem]{Conjecture}

\theoremstyle{definition}

\newtheorem{Remark}[Theorem]{Remark}

\numberwithin{equation}{section}

\newcounter{nootje}
\renewcommand\check[1]
  {\marginpar{\tiny\begin{minipage}{20mm}\begin{flushleft}\thenootje : #1\end{flushleft}\end{minipage}}\addtocounter{nootje}{1}}
\begin{document}

\title[Rational points on curves of genus 2]%
      {On the average number of rational points \\ on curves of genus 2}

\author{Michael Stoll}
\address{Mathematisches Institut,
         Universit\"at Bayreuth,
         95440 Bayreuth, Germany.}
\email{Michael.Stoll@uni-bayreuth.de}
\date{\today}

\maketitle


\section{Introduction}

For $N > 0$, let $\CC_N$ denote the set of all genus~$2$ curves
\[ C : y^2 = F(x,z) = f_6\,x^6 + f_5\,x^5 z + \dots + f_1\,x z^5 + f_0\,z^6 \]
with integral coefficients $f_j$ such that $|f_j| \le N$ for all~$j$.
($C$ is considered in the weighted projective plane with weights~$1$
for $x$ and~$z$ and weight~$3$ for~$y$.)

In this note, we sketch heuristic arguments that lead to
the following conjectures.

\begin{Conjecture} \label{Conj}
  There is a constant $\gamma > 0$ such that
  \[ \frac{\sum_{C \in \CC_N}\# C(\Q)}{\# \CC_N} \sim \frac{\gamma}{\sqrt{N}}
     \,. 
  \]
  In particular, the density of genus~$2$ curves with a rational point
  is zero.
\end{Conjecture}

The second part of this conjecture is analogous to Conjecture~2.2~(i)
in~\cite{PoonenVoloch}, which considers hypersurfaces in~$\BP^n$.

If $C$ is a curve of genus~2 as above and $P = (a : y : b)$ is a rational point
on~$C$ (i.e., we have $F(a,b) = y^2$ with $a, b$ coprime integers), then
we denote by $H(P)$ the height $H(a:b) = \max\{|a|,|b|\}$ of its $x$-coordinate.

\begin{Conjecture} \label{ConjBound}
  Let $\eps > 0$. Then there is a constant $B_\eps$ and a Zariski open
  subset~$U_\eps$ of the `coefficient space' $\BA^7$ such that for all
  $C \in \CC_N \cap U_\eps$ and all rational points $P$ on~$C$, we have
  \[ H(P) \le B_\eps N^{13/2 + \eps} \,. \]
\end{Conjecture}

The reason for restricting to~$U_\eps$ is that one should expect infinite
families of curves with larger points (at least over sufficiently large
number fields). In general, we still expect the following to hold.

\begin{Conjecture} \label{ConjBoundGen}
  There are constants~$\kappa$ and~$B$ such that every rational point~$P$
  on any curve $C \in \CC_N$ satisfies $H(P) \le B N^\kappa$.
\end{Conjecture}

If we restrict to quadratic twists of a fixed curve, then the ABC~Conjecture
implies such a bound with $\kappa = 1/2$, see~\cite{Granville}.

Note that Conjecture~\ref{ConjBoundGen} says in particular that the height
of a point on~$C$ is
{\em polynomially bounded} by the height of~$C$. If a statement like the
above could be proved for some explicit~$\kappa$ and~$B$,
then this would immediately imply that there is a polynomial
time algorithm that determines the set of rational points on a given
curve~$C$ of genus~2. More precisely, it would be polynomial time in~$N$
(and not in the input length, which is roughly~$\log N$). If we assume
that the Mordell-Weil group of the Jacobian~$J$ of~$C$ is known, then we
obtain a very efficient algorithm, since we only have to check all points
in~$J(\Q)$ of logarithmic height $\ll \log N$.

Similar statements can be formulated for other families of curves.

We also present the conjecture below, which is based on observation of
experimental data, and not on our heuristic arguments.

\begin{Conjecture} \label{ConjNumber}
  There is a constant~$B$ such that for any curve $C \in \CC_N$, the
  number of rational points on~$C$ satisfies
  \[ \#C(\Q) \le B \log (2N + 1) \,. \]
\end{Conjecture}

Caporaso, Harris, and Mazur~\cite{CHM} show that the weak form of Lang's
conjecture on rational points on varieties of general type (namely, that
they are not Zariski dense) would imply that there is a uniform bound
on~$\#C(\Q)$, independent of~$N$. So our conjecture here can be considered
as a weaker form of this consequence of Lang's conjecture.

\subsection*{Acknowledgments}

I thank Noam Elkies for providing me with his wonderful ternary sextics.
I also wish to thank Noam Elkies, Bjorn Poonen and Samir Siksek for some
helpful comments on earlier versions of this text.


\section{The Heuristic}

We first need an estimate for the fraction of curves of the form
$y^2 = F(x,z)$ in a $(1,3,1)$-weighted projective plane, with $F$
a sextic form with integral coefficients bounded by~$N$ in absolute
value, that are singular (and so are not of genus~$2$). The corresponding
forms~$F$ have a repeated irreducible factor. The largest contribution
to the set $\CD_N$ of singular curves comes from polynomials with a
repeated linear factor; they are of the form
\[ F(x, z) = (ax + bz)^2 G(x, z) \]
with $\deg G = 4$, with coefficients such that $F(x, z)$ has coefficients
bounded by~$N$. For fixed $(a:b)$, we denote by $H(a : b) = \max\{|a|,|b|\}$
the usual height in~$\BP^1$; then this number is bounded by (roughly)
$(2N+1)^5/H(a:b)^{10}$, leading to $\#\CD_N = O(N^5)$. Hence
$\#\CC_N/(2N+1)^7 = 1 - O(N^{-2})$. See Section~\ref{Sbad} below for
details.

We try to estimate the average number of rational points on
curves in~$\CC_N$ with given $x$-coordinate $(a : b) \in \BP^1(\Q)$.
Denote this number by $\E_{(a:b)}(N)$.
In the simplest
case, $(a : b) = (1 : 0)$ (or $(0 : 1)$, which leads to the same
computation). For a given curve (identified with the sextic form~$F$)
to have such a rational point, its coefficients have to satisfy
\[ f_6 = u^2 \qquad\text{for some $u \in \Z_{\ge 0}$.} \]
If $u = 0$, we have one point, for $u > 0$, we have two.
The total number of such points on (not necessarily nonsingular)
curves $y^2 = F(x,z)$ is then
\[ (2\lfloor\sqrt{N}\rfloor + 1)(2N + 1)^6 \,. \]
The number of all polynomials is $(2N + 1)^7$, and if we neglect those
that are not squarefree (which is allowed, see above), we obtain
for the average number of points at infinity
\[ \E_{(1:0)}(N) = \frac{2\lfloor\sqrt{N}\rfloor + 1}{2N+1}
             \sim \frac{1}{\sqrt{N}} \,. \]
For $(a:b) \neq (1:0), (0:1)$, we claim that similarly (see Cor.~\ref{Corab})
\begin{equation} \label{Claimab}
  \E_{(a:b)}(N) \sim \frac{\gamma(a:b)}{\sqrt{N}} 
\end{equation}
with, for $0 < a < b$,
\[ \gamma(a:b) = \frac{1}{b^3}\,\phi\Bigl(\frac{a}{b}\Bigr) \,, \]
where, for $t > 0$,
\[ \phi(t) = \frac{1}{3 \cdot 5 \cdot 7 \cdot 9 \cdot 11 \cdot 13\,t^{21}}
             \sum_{\eps_0,\dots,\eps_6 \in \{\pm 1\}}
              \eps_0 \eps_1 \cdots \eps_6
              \max\{\eps_0 + \eps_1 t + \dots + \eps_6 t^6, 0\}^{13/2} \,.
\]
In general, we have 
$\gamma(a : b) = \gamma\bigl(\min\{|a|,|b|\} : \max\{|a|,|b|\}\bigr)$.

\begin{figure}[ht] 
\begin{center}
  \Gr{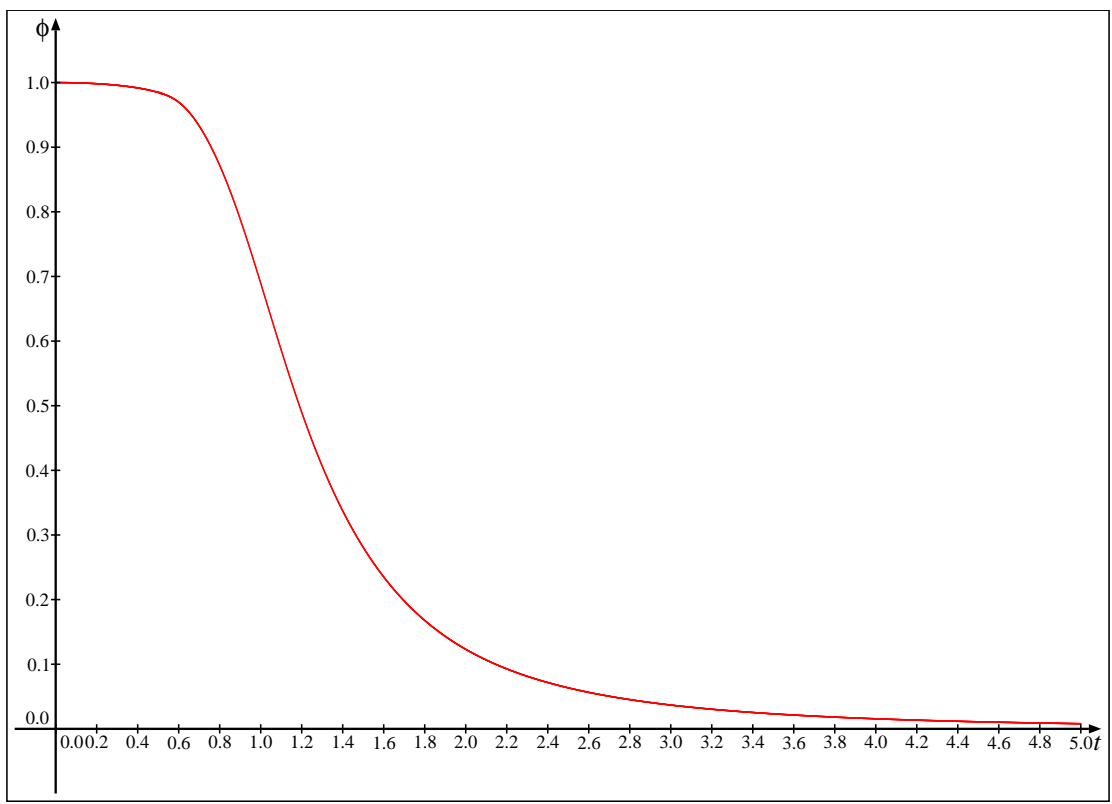}{\textwidth}
\end{center}
\caption{The function $\phi$. For $0 \le t \le 0.5$, the power series was used,
         for $0.5 \le t \le 2$ the sum, and for $t \ge 2$ the functional equation.}%
        \label{FigPhi}
\end{figure}

Note that for $t \to 0$,
\begin{align*}
   \phi(t) &= 1 - \frac{1}{2^3 \cdot 3}\,t^2
                - \frac{19}{2^7 \cdot 3}\,t^4 
                - \frac{217}{2^{10} \cdot 3}\,t^6 
                - \frac{9583}{2^{15} \cdot 3}\,t^8 
                - \frac{40125}{2^{18}}\,t^{10} 
                + O(t^{12}) \,,
\end{align*}
so that we can extend $\phi$ to all of~$\R$ by setting $\phi(0) = 1$
and $\phi(t) = \phi(|t|)$. The power series expansion is obtained by
noting that for $|t| \le 1/2$, we have
\[ \eps_0 + \eps_1 t + \dots + \eps_6 t^6 \ge 0 
   \quad\iff\quad \eps_0 = +1 \,.
\]
The radius of convergence of the series is given by the positive root
$\rho \approx 0.504138$ of $1 - t - t^2 - \dots - t^6$. 
We have the functional equation (for $t \neq 0$)
\[ \phi\Bigl(\frac{1}{t}\Bigr) = t^3\,\phi(t) \,. \]
Furthermore, $\phi(t)$ is decreasing for $t \ge 0$. This implies that
\[ \frac{\phi(1)}{H(a:b)^3} \le \gamma(a:b) \le \frac{1}{H(a:b)^3} \,. \]
Note that
\[ \phi(1) = \frac{7^{13/2} - 7\cdot 5^{13/2} + 21\cdot 3^{13/2} - 35}{135135} 
           \approx 0.689540287634369059265 \,. 
\]
See Figure~\ref{FigPhi} for a graph of~$\phi$.

We postpone the proof of the claim~\eqref{Claimab} to Section~\ref{Spfab}.

Summing the terms for $H(a:b) \le H$, we obtain, denoting by
$\E_{\le H}(N)$ the average number of rational points of height $\le H$
(where the height of a rational point is the usual naive height
$H(a:b) = \max\{|a|, |b|\}$ of its $x$-coordinate $(a:b)$):
\[ \E_{\le H}(N) \sim \frac{\gamma_H}{\sqrt{N}} \qquad
     \text{as $N \to \infty$, uniformly for $H \ll N^{6/5-\eps}$}, 
\]
where 
\[ \gamma_H = \sum_{H(a:b) \le H} \gamma(a:b) \,. \]
See Cor.~\ref{CorH} in Section~\ref{Spfab} below.

We obtain Conjecture~\ref{Conj} by letting $H \to \infty$, with
\[ \gamma = \lim_{H \to \infty} \gamma_H
          = \sum_{(a:b) \in \BP^1(\Q)} \gamma(a:b) \,.
\]
We denote by $\E(N)$ the average number of rational points on curves in~$\CC_N$.
Note that we can at least prove the following (which is, however, the less
interesting inequality).

\begin{Proposition}
  We have 
  \[ \liminf_{N \to \infty} \sqrt{N}\,\E(N) \ge \gamma \,. \]
\end{Proposition}

\begin{proof}
  Given $\eps > 0$, fix $H$ such that $\gamma_H > \gamma - \eps$.
  We then have 
  \[ \sqrt{N}\,\E(N) \ge \sqrt{N}\,\E_H(N) > \gamma_H - \eps > \gamma - 2 \eps
     \qquad\text{for $N$ sufficiently large.}
  \]
\end{proof}

In order to prove Conjecture~\ref{Conj}, one would need a reasonably
good estimate for the number of very large points. This is most likely
a very hard problem.

Let us look a bit closer at the value of~$\gamma$. We have
\[ \gamma = 4 \sum_{b=1}^\infty \sump_{0 \le a \le b, a \perp b}
                   \frac{1}{b^3} \phi\Bigl(\frac{a}{b}\Bigr)
          = \frac{4}{\zeta(3)} \sum_{H=1}^\infty 
              \frac{1}{H^3} \sump_{0 \le a \le H} \phi\Bigl(\frac{a}{H}\Bigr) 
            \,.
\]
Here, $\sum'$ denotes the sum with first and last terms counted half.

By the Euler-Maclaurin summation formula,
\[ \sump_{0 \le a \le H} \phi\Bigl(\frac{a}{H}\Bigr)
     = H \int_0^1 \phi(t)\,dt + \frac{1}{12 H} \phi'(1)
         - \frac{1}{720 H^3} \phi'''(1) + O\Bigl(\frac{1}{H^5}\Bigr) \,.
\]
So we obtain
\[ \gamma = 4\Bigl(\frac{\zeta(2)}{\zeta(3)} \int_0^1 \phi(t)\,dt 
                     + \frac{\phi'(1) \zeta(4)}{12 \zeta(3)}
                     - \frac{\phi'''(1) \zeta(6)}{240 \zeta(3)} 
                     + R \Bigr)
    \,,
\]
with a small error~$R$.

For more precise numerical estimates, we compute the first few terms
in the series over~$H$ to some precision and estimate the tail of the
series by the formula above. Note that the derivatives of~$\phi$ at~$t = 1$
can be computed explicitly. We find
\[ \gamma \approx 4.79991101188445188 \,. \]
Here is a table with experimental data obtained from all curves of size
$N \le 10$. For $N \le 3$, the number of points should be accurate; for
$4 \le N \le 10$, we counted all points of height up to $2^{14} - 1$,
so the numbers given are lower bounds. However, the difference is likely
to be so small that it does not affect the leading few digits. See
Section~\ref{S:Data} for the source of these data.
\begin{center}
\begin{tabular}{|l||c|c|c|c|c|c|c|c|c|c|}
  \hline 
  size of curves $\le N${\Large\strut}
    & 1 & 2 & 3 & 4 & 5 & 6 & 7 & 8 & 9 & 10 \\\hline
  avg.~$\#C(\Q)${\Large\strut}
    & 3.94 & 2.70 & 2.19 & 2.42 & 2.08 & 1.84 & 1.66 & 1.52 & 1.65 & 1.53 \\\hline
  (avg.~$\#C(\Q)$)$\sqrt{N}${\Large\strut}
    & 3.94 & 3.82 & 3.79 & 4.84 & 4.66 & 4.50 & 4.40 & 4.31 & 4.94 & 4.83 \\\hline
\end{tabular}
\end{center}
We observe values reasonably close to the expected asymptotic value
$\gamma \approx 4.800$. When $N$ is a square, the average number of points
jumps up because of the additional possibilities for points at $x = 0$
or $x = \infty$ (leading or trailing coefficient equal to~$N$).

From the above, we also get an estimate for $\gamma - \gamma_H$:
\[ \gamma - \gamma_H = 4 \sum_{b>H} \sum_{0<a<b, a \perp b} 
                         \frac{1}{b^3} \phi\Bigl(\frac{a}{b}\Bigr)
     \approx \frac{4}{\zeta(2) H}  \int_0^1 \phi(t)\,dt
     \approx 2.28253672259903912\,\frac{1}{H} \,.
\]


\section{Proof of the asymptotics for fixed $(a:b)$} \label{Spfab}

The total number of rational points with $x$-coordinate $(a : b) \in \BP^1(\Q)$
on curves in~$\CC_N \cup \CD_N$ is the number of integral solutions
$(f_0, f_1, \dots, f_6, y)$ of the equation
\[ f_6 a^6 + f_5 a^5 b + f_4 a^4 b^2 + f_3 a^3 b^3 + f_2 a^2 b^4 
    + f_1 a b^5 + f_0 b^6 = y^2 \,,
\]
subject to the inequalities $- N \le f_j \le N$ for $j = 0, 1, \dots, 6$.
If we fix~$y$, then the solutions correspond to the lattice points in
the intersection of the cube $[-N,N]^7$ with the hyperplane given by the
equation above. For $y = 0$, the intersection of $\Z^7$ with the hyperplane,
which we will denote $L_{(a:b)}$, is spanned by the vectors
\[ (-a,b,0,0,0,0,0), (0,-a,b,0,0,0,0), \dots, (0,0,0,0,-a,b,0), 
   (0,0,0,0,0,-a,b) \,.
\]
We can define the lattice spanned by these vectors for any $(a:b) \in \BP^1(\R)$.
These lattices (considered up to scaling) make up
the image of the obvious map from~$\BP^1(\R)$ into the moduli space
of $6$-dimensional lattices; this image is compact since $\BP^1(\R)$ is.
This implies that all invariants of our lattices (like for example the covering radius)
can be estimated above and below by a constant times
a suitable power of the typical length $H(a:b)$ associated to the lattice.
For some of these invariants, we give explicit bounds below.

The Gram matrix of the vectors above is tridiagonal:
\[ \begin{pmatrix}
     a^2+b^2 &     -ab &       0 & \cdots &       0 \\
         -ab & a^2+b^2 &     -ab & \cdots &       0 \\
           0 &     -ab & a^2+b^2 & \cdots &       0 \\
      \vdots &  \vdots &  \vdots & \ddots &  \vdots \\
           0 &       0 &       0 & \cdots & a^2+b^2
   \end{pmatrix}
\]
(From this matrix, one can again see that the lattice has a nearly orthogonal
basis consisting of vectors of equal length.)
The covolume of the lattice (in six-dimensional volume in~$\R^7$) is 
\[ \Delta_{(a:b)} = \sqrt{a^{12} + a^{10} b^2 + \dots + b^{12}} \,; \]
we have
\[ H(a:b)^6 \le \Delta_{(a:b)} \le \sqrt{7} H(a:b)^6. \] 
The diameter of the fundamental parallelotope spanned by these vectors is 
\[ \delta_{(a:b)} = \sqrt{a^2 + 5(|a|+|b|)^2 + b^2}
                  = \sqrt{6a^2 + 10|ab| + 6b^2} 
                  \le \sqrt{22}\,H(a:b) \,.
\]
Let 
\[ \bfa_{(a:b)} = (b^6, a b^5, a^2 b^4, \dots, a^6) \in \R^7 \]
and
\[ \bfe_{(a:b)} = \frac{1}{\Delta_{(a:b)}^2}\,\bfa_{(a:b)} \]
(note that $\bfa_{(a:b)} \cdot \bfe_{(a:b)} = 1$);
then the number of points we want to count is
\[ \sum_{y \in \Z}
     \#\bigl(\Z^7 \cap [-N,N]^7 
              \cap (L_{(a:b)} + y^2 \bfe_{(a:b)})\bigr) 
     = \sum_{y \in \Z} \# S(y) \,,
\]
where we define $S(y)$ to be the set under the `$\#$' sign in the first sum. Let 
$V_{(a:b)} \subset L_{(a:b)}$ be the Voronoi cell of the lattice
$\Z^7 \cap L_{(a:b)}$; in particular it has volume~$\Delta_{(a:b)}$, and its 
translates by lattice points tessellate~$L_{(a:b)}$.
We consider $S(y) + V_{(a:b)}$. The (6-dimensional) volume of this set 
is $\#S(y) \Delta_{(a:b)}$. We use $B_r(x)$ to denote the closed
ball of radius~$r$ with center~$x$. Write
\[ W_{(a:b)}(t,\delta)
     = \begin{cases}
         \bigl\{x \in L_{(a:b)} + t \bfe_{(a:b)}
                 : B_{-\delta(x)} \subset [-1,1]^7\bigr\} \,,
           & \text{if $\delta \le 0$,} \\
         \bigl\{x \in L_{(a:b)} + t \bfe_{(a:b)}
                 : B_\delta(x) \cap [-1,1]^7 \neq \emptyset\bigr\} \,,
           & \text{if $\delta \ge 0$.}
       \end{cases}
\]
In particular, $W_{(a:b)}(t,0) = (L_{(a:b)} + t \bfe_{(a:b)}) \cap [-1,1]^7$.

There is a constant $c_0 > 0$ such that the covering radius of the lattice
$\Z^7 \cap L_{(a:b)}$ is bounded by $c_0 H(a:b)$ (see the remark above).
Writing $H = H(a:b)$ in the following, we obtain
\[ N \cdot W_{(a:b)}\Bigl(\frac{y^2}{N}, -\frac{c_0 H}{N}\Bigr)
    \subset S(y) + V_{(a:b)} 
    \subset N \cdot W_{(a:b)}\Bigl(\frac{y^2}{N}, \frac{c_0 H}{N}\Bigr)
   \,.
\]
Since
\[ \vol_6 W_{(a:b)}(t, \delta)
      = \vol_6 W_{(a:b)}(t, 0) + O(\delta) + O(\delta^6)
\]
(with $O$-constants independent of~$(a:b)$), we obtain
\[ \#S(y) \Delta_{(a:b)}
     = N^6\,\vol_6 W_{(a:b)}\Bigl(\frac{y^2}{N}, 0\Bigr) + O(H N^5) + O(H^6) \,.
\]
Therefore, using that $\Delta_{(a:b)} \asymp H^6$ and writing
$f_{(a:b)}(t) = \vol_6 W_{(a:b)}(t, 0)$,
\[ \#S(y) = \frac{N^6}{\Delta_{(a:b)}} f_{(a:b)}\Bigl(\frac{y^2}{N}\Bigr)
             + O(H^{-5} N^5) + O(1) \,.
\]
Let $S(a :b)$ be the set $\cup_{y \in \Z} S(y)$ of all relevant lattice
points. Then
\begin{align*}
  \#S(a:b) &= \sum_{y \in \Z} \#S(y) \\
           &= \frac{N^6}{\Delta_{(a:b)}}
              \sum_{y \in \Z} f_{(a:b)}\Bigl(\frac{y^2}{N}\Bigr)
                  + O(H^{-2} N^{11/2}) + O(H^3 N^{1/2}) \\
           &= \frac{2 N^6}{\Delta_{(a:b)}}
              \int_0^\infty f_{(a:b)}\Bigl(\frac{y^2}{N}\Bigr)\,dy \,
                  + O(H^{-6} N^6) + O(H^{-2} N^{11/2}) + O(H^3 N^{1/2}) \,.
\end{align*}
(Note that $y = O(\sqrt{N \Delta_{(a:b)}}) = O(H^3 \sqrt{N})$ and that
$f_{(a:b)}(t)$ is decreasing for $t \ge 0$, with $f(0) = O(1)$.) Substituting 
$t = y^2/N$, this gives
\begin{align*}
  \#S(a:b) &= \frac{N^{13/2}}{\Delta_{(a:b)}}
               \int_0^\infty f_{(a:b)}(t)\,\frac{dt}{\sqrt{t}}
                 + O(H^{-6} N^6) + O(H^{-2} N^{11/2}) + O(H^3 N^{1/2}) \\
           &= N^{13/2}
               \int_{[-1,1]^7} (\bfa_{(a:b)} \cdot \bfx)_+^{-1/2}\,d\bfx \\
           &\qquad\qquad{}+ O(H^{-6} N^6) + O(H^{-2} N^{11/2}) + O(H^3 N^{1/2}) 
           \,.
\end{align*}
Here $x_+^{-1/2}$ is zero when $x \le 0$ and $x^{-1/2}$ when $x > 0$.
More generally, for $x, r \in \R$, we let $x_+^r$ denote $0$ when $x \le 0$
and $x^r$ when $x > 0$.

\begin{Lemma} \label{Lemma:Int}
  We have, for $ab \neq 0$,
  \begin{align*}
    \int_{[-1,1]^7} &(\bfa_{(a:b)} \cdot \bfx)_+^{-1/2}\,d\bfx \\
      &= \frac{2^7}{135135 \, |a b|^{21}}\!
        \sum_{\eps_0,\dots,\eps_6 = \pm 1}\! \eps_0 \cdots \eps_6
          (\eps_0 |b^6| + \eps_1 |a b^5| + \dots + \eps_6 |a^6|)_+^{13/2} .
  \end{align*}
\end{Lemma}

Note that this is $2^7$ times
\[ \gamma(a:b) = \frac{1}{|b|^3} \phi\Bigl(\frac{|a|}{|b|}\Bigr) \]
in the notation introduced in the previous section.

\begin{proof}
  Let $a_1, \dots, a_m > 0$ be real numbers, $r > -1$, $c \in \R$.
  Write $\bfa = (a_1, \dots, a_m)$.
  We prove the more general statement
  \begin{align*}
    \int_{[-1,1]^m} &(\bfa \cdot \bfx + c)_+^r\,d\mathbf{x} \\
       &\hspace*{-3mm}{}= \frac{1}{a_1 \cdots a_m\,(r+1) \cdots (r+m)}
          \sum_{\eps_1, \dots, \eps_m = \pm 1} \eps_1 \cdots \eps_m
                           (\eps_1 a_1 + \dots + \eps_m a_m + c)_+^{r+m} \,.
  \end{align*}
  We proceed by induction. When $m = 1$, we have
  \[ \int_{-1}^1 (a_1 x_1 + c)_+^r\,dx_1
       = \frac{1}{a_1\,(r+1)} 
            \bigl((a_1 + c)_+^{r+1} - (-a_1 + c)_+^{r+1}\bigr) \,,
  \]
  as can be checked by considering the cases $-c \le -a_1$, 
  $-a_1 \le -c \le a_1$, and $a_1 \le -c$ separately.
  
  For the inductive step, we assume the statement to be true for
  $a_1, \dots, a_m$ and~$r$, and prove it for $a_1, \dots, a_m, a_{m+1}$.
  Let $\bfa' = (a_1, \dots, a_m)$ and
  $\bfa = (a_1, \dots, a_{m+1})$, and use similar notation for
  vectors $\bfx$, $\bfx'$. Then
  \begin{align*}
    &\int_{[-1,1]^{m+1}} (\bfa \cdot \bfx + c)_+^r\,d\bfx \\
      &= \int_{-1}^1 \int_{[-1,1]^m}
               \bigl(\bfa' \cdot \bfx' + a_{m+1} x_{m+1} + c\bigr)_+^r
               \,d\bfx'\,dx_{m+1} \\
      &= \int_{-1}^1 \frac{1}{a_1 \cdots a_m\,(r+1) \cdots (r+m)} \times{} \\
      &\quad \sum_{\eps_1, \dots, \eps_m = \pm 1} \eps_1 \cdots \eps_m
    \bigl(\eps_1 a_1 + \dots + \eps_m a_m + x_{m+1} a_{m+1} + c\bigr)_+^{r+m}
                     \,dx_{m+1} \\
      &= \frac{1}{a_1 \cdots a_m\,(r+1) \cdots (r+m)} \times{} \\
      &\quad \sum_{\eps_1, \dots, \eps_m = \pm 1} \eps_1 \cdots \eps_m
        \int_{-1}^1
      \bigl(\eps_1 a_1 + \dots + \eps_m a_m + x_{m+1} a_{m+1} + c\bigr)_+^{r+m}
                     \,dx_{m+1} \\
      &= \frac{1}{a_1 \cdots a_m\,(r+1) \cdots (r+m)} \times{} \\
      &\quad \sum_{\eps_1, \dots, \eps_m = \pm 1} \eps_1 \cdots \eps_m
        \frac{1}{a_{m+1}\,(r+m+1)}
         \sum_{\eps_{m+1} = \pm 1} 
          \bigl(\eps_1 a_1 + \dots + \eps_{m+1} a_{m+1}
                         + c\bigr)_+^{r+m+1}
  \end{align*}
  by the case $m = 1$.
  
  To finish the proof of the lemma, note that we can take $a, b > 0$.
  We then apply the claim with $\bfa = \bfa_{(a:b)}$, $r = -1/2$, and $c = 0$.
\end{proof}

\begin{Corollary} \label{Corab}
  With $H = H(a:b)$,
  \[ \E_{(a:b)}(N) = \frac{\gamma(a:b)}{\sqrt{N}}
                        + O\bigl(H^{-6} N^{-1}\bigr)
                        + O\bigl(H^{-2} N^{-3/2}\bigr)
                        + O\bigl(H^3 N^{-13/2}\bigr) \,.
  \]
  In particular, we have
  \[ \sqrt{N}\,\E_{(a:b)}(N) \To \gamma(a:b) \]
  as $N \to \infty$, uniformly for $(a:b)$ such that $H(a:b) \ll N^{2-\eps}$.
\end{Corollary}

\begin{proof}
  First note that $\E_{(a:b)}(N) = \#S'(a:b)/\#\CC_N$, where 
  $S'(a:b)$ only lists the points in $S(a:b)$ on curves that are smooth.
  We have 
  \[ \#\CC_N = (2N+1)^7 - \#\CD_N = (2N+1)^7 + O(N^5)
             = (2N)^7\bigl(1 + O(N^{-1})\bigr)
  \]
  and
  \[ \#S'(a:b) = \#S(a:b) 
                  + O\bigl(H^{-10} N^5\bigr)
                  + O\bigl(H^{-1} N^{9/2}\bigr) 
                  + O\bigl(H^3 N^{1/2}\bigr) \,. \]
  See Section~\ref{Sbad} below.
  This implies that
  \[ \E_{(a:b)}(N)
        = \frac{\#S(a:b)}{(2N)^7}\bigl(1 + O(N^{-1})\bigr) 
           + O\bigl(H^{-10} N^{-2}\bigr)
           + O\bigl(H^{-1} N^{-5/2}\bigr)
           + O\bigl(H^3 N^{-13/2}\bigr) \,. 
  \]
  By Lemma~\ref{Lemma:Int}, the definition of $\gamma(a:b)$, and the discussion
  preceding the lemma, we have (using $\gamma(a:b) \asymp H^{-3}$)
  \[ \frac{\#S(a:b)}{(2N)^7}
       = \frac{\gamma(a:b)}{\sqrt{N}}
           + O\bigl(H^{-6} N^{-1}\bigr)
           + O\bigl(H^{-2} N^{-3/2}\bigr)
           + O\bigl(H^3 N^{-13/2}\bigr) \,.
  \]
  The result is obtained by combining these results, after eliminating
  redundant terms.
\end{proof}

\begin{Corollary} \label{CorH}
  \[ \E_{\le H}(N) = \frac{\gamma_H}{\sqrt{N}}
                        + O\bigl(N^{-1}\bigr)
                        + O\bigl((\log H) N^{-3/2}\bigr)
                        + O\bigl(H^5 N^{-13/2}\bigr) \,.
  \]
  In particular, we have
  \[ \sqrt{N}\,\E_{\le H(N)}(N) \To \gamma_{H(N)} \]
  as $N \to \infty$ if $H(N) \ll N^{6/5-\eps}$, and
  \[ \E_{\le H(N)}(N) \To 0 \]
  as $N \to \infty$ if $H(N) \ll N^{13/10-\eps}$.
\end{Corollary}

\begin{proof}
  Sum the estimates in the previous corollary.
\end{proof}

It should be possible to extend the range beyond $H \ll N^{6/5 - \eps}$ if one
uses more sophisticated methods from analytic number theory. (In fact,
Stephan Baier~\cite{Baier} has obtained an exponent of $7/5-\eps$.)
It would be interesting to see how far one can get.


\section{Counting Bad Curves and Points} \label{Sbad}

In this section, we will bound the number $\#\CD_N$ of non-smooth curves
and the total number of points of height $\le H$ on them. Recall the
following.

\begin{Lemma} \label{Lemma:GenBound}
  Let $\Lambda \subset \R^n$ be a lattice of covolume~$\Delta$ and
  covering radius~$\rho$. Let $S \subset \R^n$ be a subset. Then
  \[ \#(S \cap \Lambda) \le \frac{\vol (S + B_\rho(0))}{\Delta} \,. \]
\end{Lemma}

\begin{proof}
  Let $V$ be the Voronoi cell of~$\Lambda$ (centered at zero), then
  $V \subset B_{\rho}(0)$ by definition of the covering radius, and
  $\vol V = \Delta$. It follows that
  \[ \bigcup_{x \in S \cap \Lambda} (V + x) \subset S + B_{\rho}(0)\,,
       \quad\text{and thus}\quad
     \Delta \cdot \#(S \cap \Lambda) \le \vol (S + B_{\rho}(0)) \,.
  \]
\end{proof}

To make life a bit simpler, we observe that 
$[-N,N]^7 \subset B_{\sqrt{7}N}(0)$; we will bound the number of bad curves
in the ball. This has the advantage that the intersection with any affine
subspace will be a ball again.

Note that a form $F(x,z)$ is not square-free if and only if it is divisible by
the square of a primitive form~$G$. Let $n$ be the degree of~$G$; 
assume it has coefficients $\bfal = (\alpha_n, \dots, \alpha_0)$. Then
the forms divisible by~$G^2$ correspond to lattice points in the span
of 
\[ x^{6-2n} G^2\,, x^{5-2n} z G^2\,, \dots\,, z^{6-2n} G^2\,, \] 
intersected with the ball $B_{\sqrt{7}N}(0)$. 

We can extend this to $G$ with real coefficients; then the lattices we
obtain (modulo scaling) are parametrized by the compact set $\BP^n(\R)$,
hence they all live in a compact subset of the moduli space of lattices.
Taking into account that the basis vectors have length
of order $H(\bfal)^2$, this gives the
following relations for the covolume, covering radius and minimal length
of the lattices.
\[ \Delta \asymp H(\bfal)^{14-4n}\,, \qquad \rho \asymp H(\bfal)^2\,, \qquad
   \mu \asymp H(\bfal)^2 \,.
\]
In particular, there will be no non-zero lattice point in the ball
of radius $\sqrt{7} N$ when $N > \text{const}\,H(\bfal)^2$.
By Lemma~\ref{Lemma:GenBound}, we then obtain a bound
\begin{align*}
   \#\CD_N
     &\le \sum_{n=1}^3 
           \sum_{\bfal \in \BP^n(\Q), H(\bfal) \ll \sqrt{N}}
            O\Bigl(\frac{(N + H(\bfal)^2)^{7-2n}}{H(\bfal)^{14-4n}}\Bigr) \\
     &= \sum_{n=1}^3
           \sum_{\bfal \in \BP^n(\Q), H(\bfal) \ll \sqrt{N}}
            O\Bigl(\frac{N^{7-2n}}{H(\bfal)^{14-4n}}\Bigr) \\
     &= \sum_{n=1}^3 N^{7-2n} \sum_{H \ll \sqrt{N}} O(H^{3n-14})
      = O(N^5) \,.
\end{align*}
We conclude that in fact $\#\CD_N \asymp N^5$, since we already get
$N^5$ from $G = x$.

\medskip

Now in order to count points on these bad curves, we use the same basic
idea as before. This time, we have to count lattice points in the ball
that are in a translate of the subspace of forms that are divisible
by $G(x,z)^2 (bx-az)$. We assume for now that $G(a,b) \neq 0$.
If $F(a,b) = G(a,b)^2 y^2$, then the translation
is by a vector of length $G(a,b)^2 y^2/\Delta_{(a:b)}$. 
So for the count of points with $x$-coordinate $(a:b)$, we get a bound of
\[ \sum_{|y| \ll \frac{\sqrt{N \Delta_{(a:b)}}}{|G(a,b)|}}
     O\Bigl(\frac{(N + H(\bfal)^2)^{6-2n}}{H(\bfal)^{12-4n} H(a:b)^{6-2n}}\Bigr)
    = O\Bigl(\frac{N^{\frac{13}{2}-2n}}%
                  {|G(a,b)| H(\bfal)^{12-4n} H(a:b)^{3-2n}}\Bigr) \,.
\]
Estimating $|G(a,b)| \ge 1$ trivially, we obtain for the total number
of such points the bound
\begin{align*}
  \sum_{n=1}^3 &\sum_{H(\bfal) \ll \sqrt{N}}
    O\Bigl(\frac{N^{\frac{13}{2}-2n}}{H(\bfal)^{12-4n} H(a:b)^{3-2n}}\Bigr) \\
     &= O\Bigl(\frac{N^{9/2}}{H(a:b)}\Bigr) + O\bigl(N^{5/2} H(a:b)\bigr)
        + O\bigl(N^{1/2} H(a:b)^3\bigr) \,.
\end{align*}
The middle term is redundant, since it is always dominated by one of the
others.

If $G(a,b) = 0$, then (since we can assume $G$ to be irreducible) $n = 1$,
and we have to count all forms divisible by~$G^2$. This adds a term
of order $N^5/H(a:b)^{10}$.

\begin{Remark}
  With a similar computation as above, one can show that the number of
  curves in~$\CC_N$ with reducible polynomial~$F$ is~$O(N^6)$. Therefore
  the contribution of such curves is negligible.
\end{Remark}


\section{Speculations on the Size of Points}

Recall that
\[ \gamma_H \approx \gamma - \frac{c}{H} \]
where $c \approx 2.28253672259903912$.
If we assume Conj.~\ref{Conj}, then
the calculations above suggest that the number of curves in~$\CC_N$
that have a rational point of $x$-height $> H$ is roughly
$c N^{13/2}/H$, at least as long as $H$ is not too large compared
to~$N$, see Cor.~\ref{CorH}. If we recklessly extend
this to large~$H$, this would predict that the largest rational point
on a curve from~$\CC_N$ should have height $\ll N^{13/2+\eps}$.
One has to be careful, however, as was pointed out to me by Noam Elkies,
mentioning the case of integral points on elliptic curves as an analogy.
Considering curves in short Weierstrass form $y^2 = x^3 + Ax + B$ with
$A, B \in \Z$, heuristic considerations like those presented here predict
that integral points should be of size $\ll \max\{|A|^{1/2}, |B|^{1/3}\}^{10 + \eps}$,
but there are families that reach an exponent of~$12$. See the information
given at~\cite{ElkiesWeb}. This leads to Conjecture~\ref{ConjBound}.

Regarding possible families with larger points, we consider the case
that the coefficients $f_j$ are linear forms in the coordinates $(t:u)$
of~$\BP^1$, the coordinates $x$ and~$z$ of the point we are looking for
are homogeneous polynomials of degree~$m$ (to be determined), and
$y^2 = q(t,u)^2 r(t,u)$ with $q$ of degree~$3m-1$ and $r$ of degree~$3$.
If we find a solution of
\[ q^2 r = \sum_{j=0}^6 f_j x^j z^{6-j} \]
in such polynomials, then we should obtain an infinite family of curves
with points satisfying $H(P) \gg N^m$. (Of course, we have to exclude
degenerate solutions.) To see this, multiply by $r(1,0)$ (which we can
assume to be nonzero after a suitable change of coordinates on~$\BP^1$).
The equation
\[ r(1,0) r(t,u) = w^2 \]
then has the solution $(t,u,w) = (1,0,r(1,0))$, which must be contained
in a family of solutions that is parametrized by a genus~$0$ curve
(compare \cite{DarmonGranville,Beukers}). If we plug in this parametrization,
we obtain a one-dimensional family of suitable curves with base~$\BP^1$.

There are
\[ 3m + 4 + 7 \cdot 2 + 2 \cdot (m+1) = 5m + 20 \]
unknown coefficients involved in the equation above.
On the other hand, there is an action of
$\GL_2 \times \GL_2 \times \Gm$ (given by the automorphisms of $\BP^1_{(t:u)}$,
the automorphisms of $\BP^1_{(x:z)}$ (acting on $x$, $z$, and the~$f_j$ and
leaving the value of the right hand side unchanged), and scaling of $q$
versus~$r$), which takes away 9~degrees of freedom. The relation above
leads to $6m+2$ equations, so the remaining number of degrees of freedom
should be
\[ (5m + 20) - 9 - (6m + 2) = 9 - m \,. \]
This suggests that there should be families of curves with points
such that $H(P) \gg N^9$. (We do not get better results when we take
coefficients~$f_j$ of higher degree, taking $\deg r = 4$ in case this
degree is even.) Of course, the corresponding variety may
fail to have rational points, so that we do not see these families
over~$\Q$. Or some other accidents can occur, leading to extraneous
solutions with larger~$m$.


In the following, we will ignore such special families and try to
make our `generic' conjecture more precise by using a probabilistic model.
In this model, we interpret the quantity
\[ \frac{N^6}{\Delta_{(a:b)}} \int_0^\infty f_{(a:b)}\Bigl(\frac{y^2}{N}\Bigr) \,dy
    = 2^6 \gamma(a : b) N^{13/2}
\]
that gives rise to the main term in the count of points $(a:\pm y:b)$
as the probability that such a point pair occurs in~$\CC_N$. The number of
pairs of points of height $> H$ should then follow a Poisson distribution
with mean 
\[ \mu_H = 2^6 (\gamma - \gamma_H) N^{13/2}
         \approx \frac{2^6 c N^{13/2}}{H} \,, 
\]
at least when $H$ is large compared to~$N$. Taking $H = \lambda N^{13/2}$,
the probability that no such point exists is then $e^{-2^6 c/\lambda}$.
Taking into account the fact that points occur in packets of eight%
\footnote{We consider all points on all curves of fixed size together.}
(change the sign of $x$ or~$y$, send $x$ to $1/x$), i.e., four point
pairs, we should correct this to $e^{-16 c/\lambda}$. For a fifty-fifty
chance of no larger points, we should take $\lambda \approx 53$, for
an 80\% chance, we take $\lambda \approx 164$. This line of argument
would lead us to expect the following.

{\em
  If $\lambda(N) \to \infty$ as $N \to \infty$, then there are only
  finitely many `generic' curves $C \in \CC_N$ of genus~2 such that $C$
  has a rational point~$P$ with $H(P) > \lambda(N) N^{13/2}$.
}

The problem with this is that it is not so clear how to make the
restriction to `generic' curves precise. There might be an infinity
of families of curves with points of height $N^k$ for a sequence of~$k$
tending to~$13/2$ from above, which could lead to problems when $\lambda(N)$
tends to infinity very slowly. Therefore, we keep on the safe side
with the given formulation of Conjecture~\ref{ConjBound}.

On the other hand, we can use similar heuristic arguments for any
given family of curves. It is reasonable to expect that the bounds we
obtain will not get arbitrarily large (in terms of the exponent of~$N$
in the height bound). This leads to Conjecture~\ref{ConjBoundGen}.

\begin{figure}[ht] 
\begin{center}
  \Gr{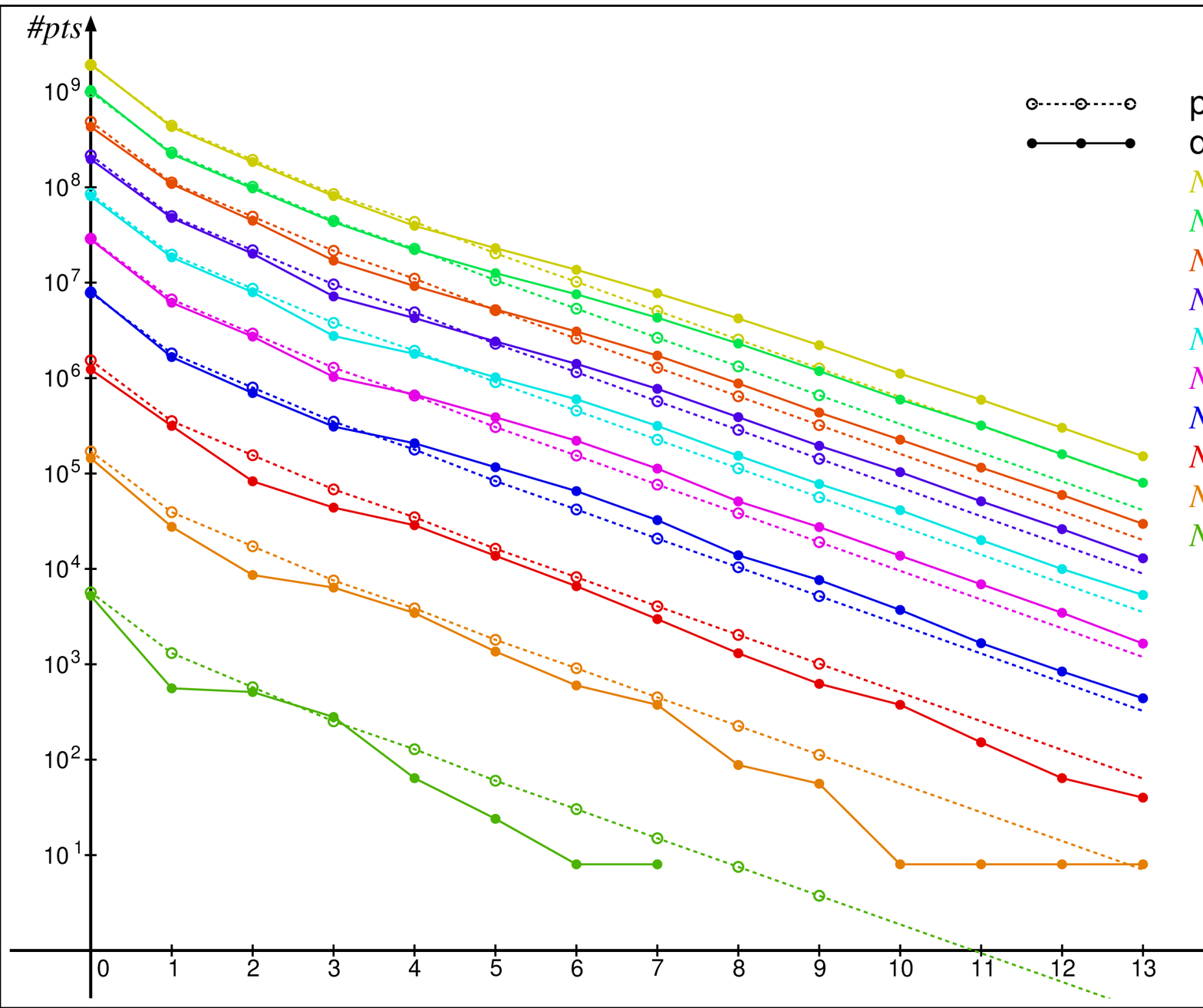}{\textwidth}
\end{center}
\caption{Expected and actual number of rational points in various height
         brackets, for $1 \le N \le 10$.}\label{FigHeights}
\end{figure}

We have checked experimentally how well the expected number of points
of height in the interval $\left[2^n, 2^{n+1}\right[$ matches the
actual number of points on curves of small size. For values of~$n$ that
are not very small, this is $2^{-(n+1)} c\,\#\CC_N/\sqrt{N}$.
Figure~\ref{FigHeights}
shows this comparison, for curves in $\CC_N$, for $1 \le N \le 10$
and $0 \le n \le 13$.
The fit is quite good, even though the range of~$N$ is certainly far
too small for the asymptotics to kick in except for very small heights.
There is an unexpected feature: starting with $N = 4$, points of larger
height seem to occur more frequently than they should.%
\footnote{Since points accumulate on singular curves, which we did not
consider here, one would perhaps rather expect a deviation in the other
direction!}
It would be
interesting to find an explanation for this phenomenon. One possibility
is that it might be related to the existence of families of curves
with systematically occurring large points. Of course,
according to our results, this can only occur for fairly large heights
when $N$ is large. See Section~\ref{S:Data} for a description of the
computations.

It is also interesting to compare the observed value of~$\lambda(N)$ such
that no rational point of height $> \lambda(N) N^{13/2}$ exists on a curve
in~$\CC_N$ with the estimates given above. For $N = 1, 2, 3$, the largest
points we found on curves in~$\CC_N$ have heights as follows.

\begin{center}
\begin{tabular}{|l||c|c|c|}
  \hline 
  size of curves{\large\strut} & $N = 1$ & $N = 2$ & $N = 3$ \\\hline
  max.~$H(P)${\large\strut}    & 145   & 10711 & 209040 \\\hline
\end{tabular}
\end{center}

We therefore find
\[ \lambda(1) \approx 145.00\,, \quad 
   \lambda(2) \approx 118.34\,, \quad
   \lambda(3) \approx 165.55\,,
\]
corresponding to probabilities (for no larger point to exist, in the
sense explained above) between 73\% and~81\%.

The record point on \quad $y^2 = x^6 - 3 x^4 - x^3 + 3 x^2 + 3$ \quad has
$x = -\frac{58189}{209040}$.

\medskip

Similar considerations for general hyperelliptic curves of genus $g \ge 2$
lead to a heuristic estimate of
$O(N^{(4g+5)/2}/H^{g-1})$ for the number of curves with a point of
height $> H$. Therefore we would expect the points to be generically of height
\[ H \ll N^{(4g+5)/(2g-2)+\eps} = N^{2 + \frac{9}{2(g-1)} + \eps} \,. \]


\section{Speculations on the Number of Points}

We can also try to extract some information of the number of points
(or point pairs) on hyperelliptic curves. Since the linear conditions
on the coefficients coming from up to seven distinct $x$-coordinates
are linearly independent, we would expect the following.

Let $R^{(m)}_N$ be the subset of~$\CC_N$ of curves that have at least
$m$ pairs of rational points (i.e., points with $m$ distinct
$x$-coordinates).
For $0 \le m \le 7$, there are constants $\gamma^{(m)} > 0$ such that
\[ \#R^{(m)}_N \sim \gamma^{(m)} \, N^{7-m/2} \,. \]


One caveat here is that the number of non-squarefree polynomials
will be in the range of these sizes if $m \ge 4$, so the conclusion is not automatic.
Indeed, the experimental data show a noticeable deviation from this
expectation already for $m \ge 3$.

Let us be more precise and try to obtain numerical values for the
$\gamma^{(m)}$. Assuming the occurrence or not of points with distinct
$x$-coordinates to be independent for all~$x \in \BP^1(\Q)$,
the generating function for the probability of
having rational points with exactly $m$ distinct $x$-coordinates
should be (assuming exact probability $\gamma(a:b)/2\sqrt{N}$ for
a point~$P \in C(\Q)$ with $x(P) = (a:b)$)
\[ G(T) = \sum_{m=0}^\infty \Prob_N\bigl(\#x(C(\Q)) = m\bigr) T^m
        = \prod_{(a:b) \in \BP^1(\Q)} 
             \Bigl(1 + \frac{\gamma(a:b)}{2\sqrt{N}} (T-1)\Bigr) \,.
\]
The numbers $\gamma^{(m)}$ should then occur as the limits as
$N \to \infty$ of the coefficients in the series
\[ \sum_{m=0}^\infty \gamma^{(m)}(N)\,T^m
     = \frac{1 - \sqrt{N}T\,G(\sqrt{N}T)}{1 - \sqrt{N}T}
     = \frac{T\,G(\sqrt{N}T) - \frac{1}{\sqrt{N}}}{T - \frac{1}{\sqrt{N}}} \,,
\]
where $\gamma^{(m)}(N) N^{-m/2}$ is an estimate for the fraction of
curves with at least $m$ point pairs.
Now, as $N \to \infty$ and coefficient-wise, this series behaves as
\[ G(\sqrt{N} T) = \prod_{(a:b)}
                    \Bigl(1 - \frac{\gamma(a:b)}{2\sqrt{N}}
                           + \frac{\gamma(a:b)}{2}\,T\Bigr) 
                 \To \prod_{(a:b)} \Bigl(1 + \frac{\gamma(a:b)}{2}\,T\Bigr) \,.
\]
So $\gamma^{(m)}$ is the degree-$m$ ``infinite elementary symmetric 
polynomial'' in the numbers $\gamma(a:b)/2$. Using $(a:b)$ of height
up to~$1000$, we find
\begin{gather*}
   \gamma^{(1)} = \frac{\gamma}{2} \approx 2.399\,, \quad
   \gamma^{(2)} \approx 2.499\,, \quad
   \gamma^{(3)} \approx 1.504\,, \quad
   \gamma^{(4)} \approx 0.591\,, \text{\quad etc.} 
\end{gather*}

\begin{figure}[ht]
\begin{center}
  \Gr{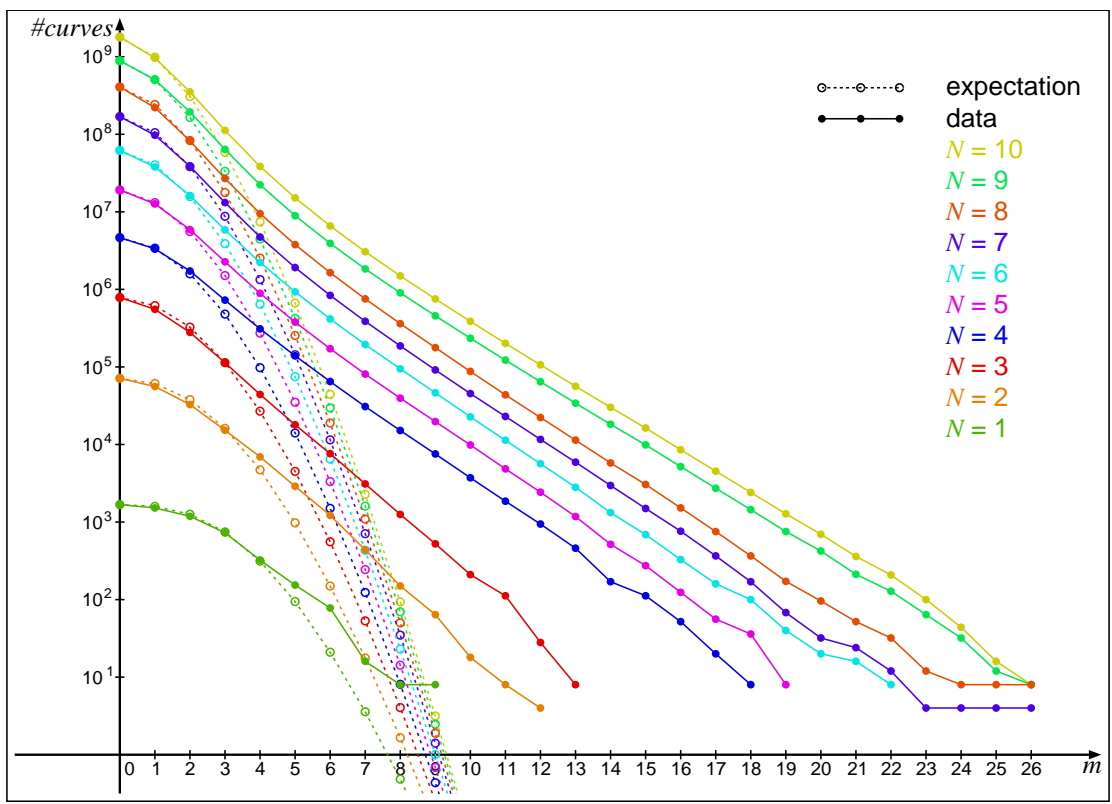}{\textwidth}
\end{center}
\caption{Expected and actual number of curves with at least $m$ pairs 
         of rational points}\label{FigNum}
\end{figure}

In Figure~\ref{FigNum}, we compare the expected values $\gamma^{(m)}(N)/N^{m/2}$
with the observed numbers. For $m \le 2$, there is good agreement, but for
$m \ge 4$, there seem to be many more curves with at least $m$ pairs of points
than predicted. Indeed, the data suggest a behavior of the form
$\alpha^{m}$ for the fraction of curves in this range, with 
$\alpha \approx 0.5$ largely independent of~$N$ (or even increasing: note
the changes in slope when $N = 4$ or $N = 9$).

This seems to indicate that as soon as there are many points, it is
much more likely that there are additional points than on average ---
the points ``conspire'' to generate more points. Maybe this is related
to another observation, which is that in examples of curves with
many rational points, the points tend to have many dependence relations
in the Mordell-Weil group. One possible explanation might be that when
there are already several points, they tend to be fairly small, so that
there are many small linear combinations of them in the Mordell-Weil group.
Such a small point in the Mordell-Weil group is represented by a pair of
points on~$C$ such that the quadratic polynomial whose roots are the
$x$-coordinates of the two points has small height. A polynomial of small
height has a good chance to split into linear factors. In this case, both
points involved are rational points on~$C$. It would be very interesting
to turn this into a precise estimate for the number~$\alpha$ that we observe.

\begin{figure}[htb]
\begin{center}
  \Gr{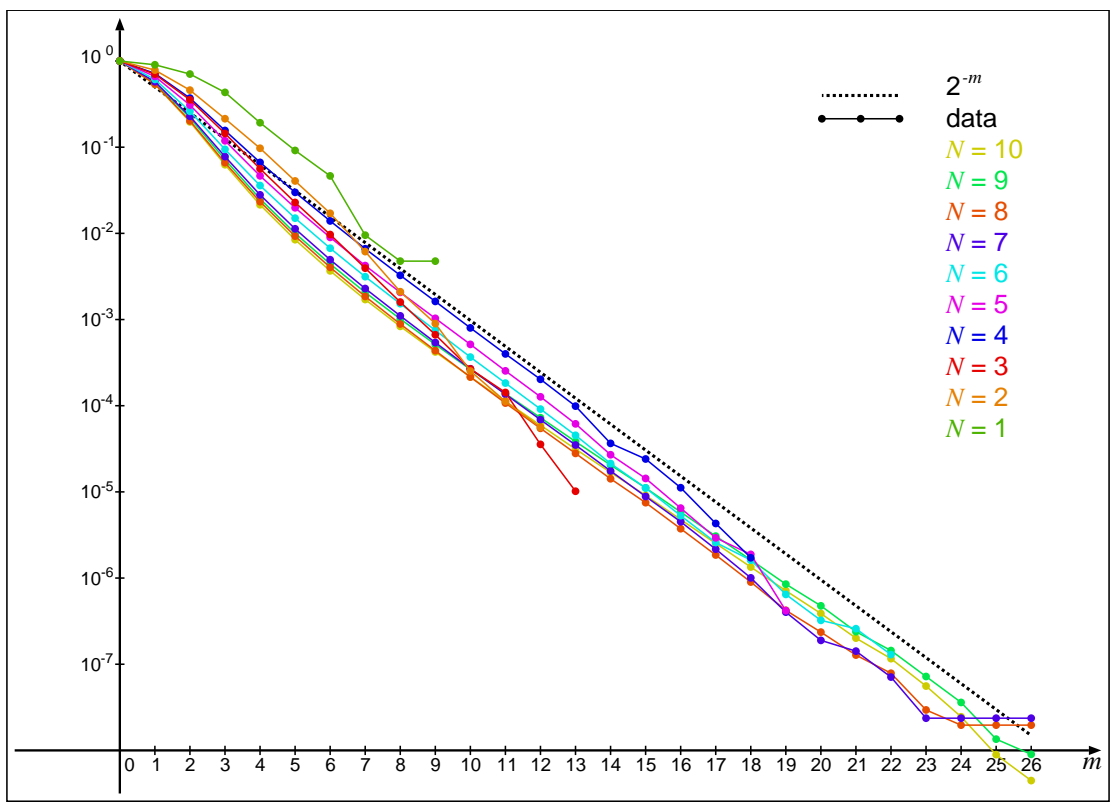}{\textwidth}
\end{center}
\caption{Proportion of curves with at least $m$ pairs of rational points}
        \label{FigNum1}
\end{figure}

In Figure~\ref{FigNum1}, we show the proportion of curves in~$\CC_N$ with
at least $m$ point pairs. It is striking how the graphs are all contained
in a narrow strip near the line (in the logarithmic scaling used in the
picture) corresponding to $m \mapsto 2^{-m}$.

If these observations extend to larger~$N$, then we should expect about
$2^{-m} (2N+1)^7$ curves in~$\CC_N$ with $m$ or more point pairs. The
largest number of point pairs on a curve in~$\CC_N$ should then be
\[ \frac{7}{\log 2}\,\log(2N + 1) + O(1) \,. \]
Conjecture~\ref{ConjNumber} gives a slightly weaker statement, replacing
the factor $7/\log 2$ by an arbitrary constant.

In order to test our conjecture, we conducted a search for curves with
many points in~$\CC_{200}$. The table in Figure~\ref{FigTable} lists the record curves
we found (curves with more point pairs than all smaller curves). On each
curve, we found all points of height up to~$2^{17} - 1 = 131\,071$
(and in some cases a few more).
The column labeled~``$F$'' lists the coefficients of one example curve.

\begin{figure}[htb] \label{FigTable}
  \[ \begin{array}{|c|c|c|c|c|} \hline
        N & F                     &  \#C(\Q) & \frac{\#C(\Q)}{\log_{10}(2N+1)}\text{\LARGE\strut}
                                      \\\hline
        1 & 1,-1, 0, 1,-1, 0, 1        &  18 & 37.73 \\
        2 & 1, 2, 0,-2, 2, 0, 1        &  24 & 34.34 \\
        3 & 1,-3, 2, 3, 0, 0, 1        &  26 & 30.77 \\
        4 & 4, 4, 0,-1,-4, 0, 1        &  36 & 37.73 \\
        5 & 4, 4, 0,-5,-4, 1, 1        &  38 & 36.49 \\
        6 & 1, 6,-1,-5, 0,-1, 1        &  44 & 39.50 \\
        7 & 4,-7, -5, 5, 1, 2, 1       &  52 & 44.21 \\
       11 & 9, 2,-11,-5, 3, 9, 9       &  56 & 41.12 \\
       13 & 9,-12,-4,13,-4, 3, 4       &  62 & 43.32 \\
       16 & 4, 1,-16,-13,16,8, 1       &  68 & 44.78 \\
       19 & 1,-18,-19,6,11,12,16       &  72 & 45.25 \\
       20 & 4, 3,20, 5,-3,-20,16       &  74 & 45.88 \\
       21 & 4, 3,19,-21,-19,14,1       &  78 & 47.75 \\
       24 & 9,24,-10,-20,2,-12,16      &  80 & 47.33 \\
       36 & 9,3,-35,5,27,-20,36        &  82 & 44.01 \\
       42 & 4,-13,23,7,-42,0,25        &  88 & 45.61 \\
       47 & 9,-21,23,-7,-47,28,16      &  98 & 49.55 \\
       54 & 9,-54,3,-2,-36,32,49       & 104 & 51.04 \\
       66 & 25,-30,-37,-46,66,34,4     & 106 & 49.91 \\
       67 & 1,-46,67,38,32,-32,4       & 114 & 53.51 \\
       70 & 49,-60,-28,-70,-9,70,49    & 118 & 54.90 \\
       72 & 1,2,63,-38,-72,36,9        & 120 & 55.52 \\
      110 & 25,-32,80,110,-105,-78,49  & 124 & 52.89 \\
      117 & 1,-26,87,83,-43,-117,64    & 126 & 53.14 \\
      125 & 49,42,-85,-125,77,69,9     & 130 & 54.17 \\
      132 & 81,-132,-16,71,76,-71,16   & 138 & 56.95 \\
      143 & 81,-120,-28,-54,143,90,9   & 140 & 56.96 \\
      184 & 1,98,-59,-184,161,46,1     & 142 & 55.32 \\
      191 & 4,-4,156,-191,-159,171,144 & 146 & 56.52 \\
      \hline
    \end{array}
  \]
  \caption{Examples of curves with many points.}
\end{figure}

\begin{figure}[htb]
\begin{center}
  \Gr{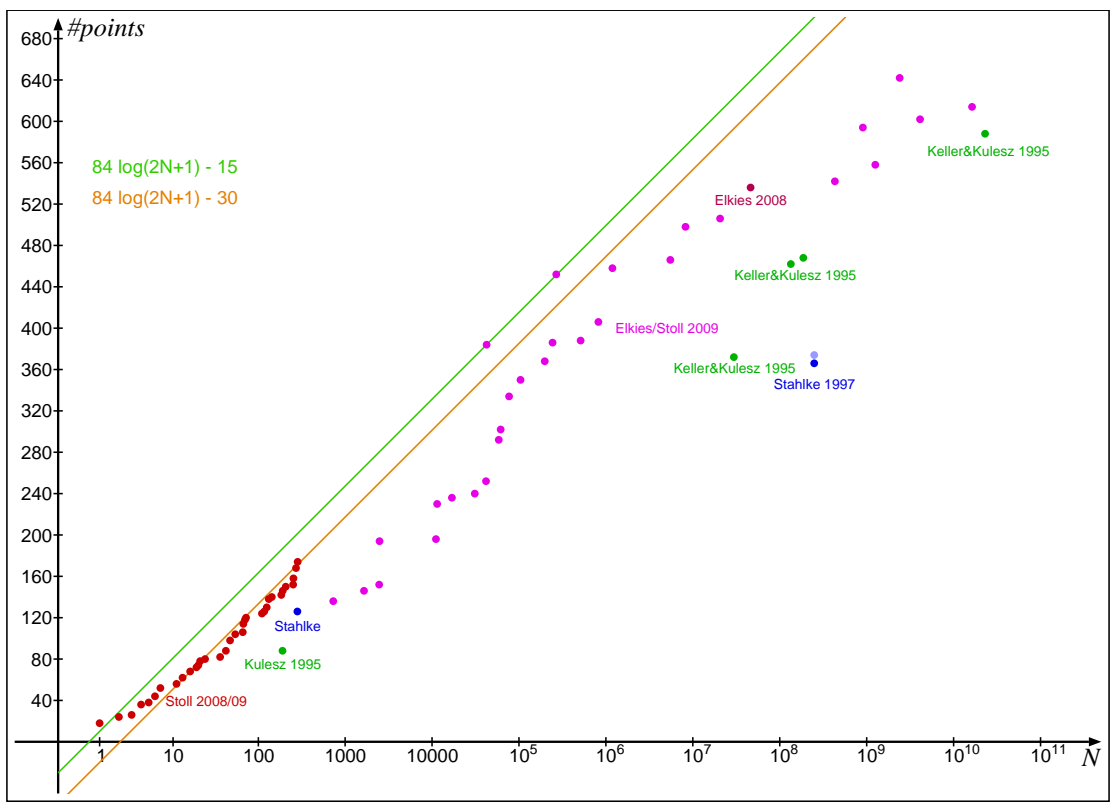}{\textwidth}
\end{center}
\caption{Curves with many points} \label{FigNum2}
\end{figure}

The constant in front of $\log(2N+1)$ that seems to fit our data best
points to a value of~$\alpha$ of about~$0.68$ in that range (corresponding to
the slope of the lines in the figure and indicating that the 
observed increase of~$\alpha$ with~$N$ persists). In Figure~\ref{FigNum2},
we have plotted $\#C(\Q)$ against $\log(2N+1)$ for the curves in the table
(and some more coming from an ongoing extended search).
In addition, we show a selection of good
curves from Elkies' families, see below, and some other previously known examples.
(``$\log$'' in the figure is the logarithm with base~$10$.) The sources
of these examples are \cite{Kulesz,KellerKulesz,Stahlke}; the curve marked
``Stahlke'' on the left was communicated to me by Colin Stahlke; it appears
in~\cite{Stoll02}, where the Mordell-Weil group of its Jacobian is determined.

One of these examples is the curve with the largest number
of point pairs found until very recently (see Keller and Kulesz~\cite{KellerKulesz}).
It has $N = 22\,999\,624\,761$ and $m = 294$. This curve has 12 automorphisms
defined over~$\Q$, and the 588~points are 49~orbits of 12~points each.
Until~2008, the record for curves with only the hyperelliptic involution
as a nontrivial automorphism was held by a curve found by Stahlke~\cite{Stahlke}
with 366 known rational points. (In fact, there are at least 8 more points,
see Section~\ref{S:Data}.)

Recently, Noam Elkies~\cite{Elkies} has constructed several K3~surfaces of the
form $y^2 = S(t,u,v)$ with a ternary sextic~$S$ such that $S$
admits a large number ($> 50$) of rational lines on which $S$ restricts
to a perfect square. Each of these therefore provides a 2-dimensional
family of genus~2 curves with more than 50 pairs of rational points.
In one of these families, he found a curve with 536~rational points.
(It is marked ``Elkies 2008'' in Figure~\ref{FigNum2}.)
In the course of a further systematic search in these families, we found
several curves with still more points, some of which even beat the Keller and Kulesz
record. The curve with the largest number of points discovered so far is
\begin{align*}
  y^2 &= 82342800 x^6 - 470135160 x^5 + 52485681 x^4 + 2396040466 x^3 \\
      & \qquad {} + 567207969 x^2 - 985905640 x + 247747600\,;
\end{align*}
it has (at least) {\bf 642} points. The $x$-coordinates of the points
with $H(P) > 10^5$ are as follows (the smaller points can easily be found
using {\tt ratpoints}, for example).
\begin{gather*}
  \frac{15121}{102391}, \frac{130190}{93793}, -\frac{141665}{55186},
  \frac{39628}{153245}, \frac{30145}{169333}, -\frac{140047}{169734},
  \frac{61203}{171017}, \frac{148451}{182305}, \frac{86648}{195399}, \\[1mm]
  -\frac{199301}{54169}, \frac{11795}{225434}, -\frac{84639}{266663},
  \frac{283567}{143436}, -\frac{291415}{171792}, -\frac{314333}{195860},
  \frac{289902}{322289}, \frac{405523}{327188}, \\[1mm]
  -\frac{342731}{523857}, \frac{24960}{630287}, -\frac{665281}{83977}, 
  -\frac{688283}{82436}, \frac{199504}{771597}, \frac{233305}{795263}, 
  -\frac{799843}{183558}, -\frac{867313}{1008993}, \\[1mm]
  \frac{1142044}{157607}, \frac{1399240}{322953},
  -\frac{1418023}{463891}, \frac{1584712}{90191}, \frac{726821}{2137953},
  \frac{2224780}{807321}, -\frac{2849969}{629081}, -\frac{3198658}{3291555}, \\[1mm]
  \frac{675911}{3302518}, -\frac{5666740}{2779443}, \frac{1526015}{5872096},
  \frac{13402625}{4101272}, \frac{12027943}{13799424}, -\frac{71658936}{86391295},
  \frac{148596731}{35675865}, \\[1mm]
  \frac{58018579}{158830656}, \frac{208346440}{37486601},
  -\frac{1455780835}{761431834}, -\frac{3898675687}{2462651894}.
\end{gather*}

The record so far for $\#C(\Q)/\log_{10}(2N+1)$ is held by the curve
\[ y^2 = 37665 x^6 - 220086 x^5 + 212355 x^4 + 268462 x^3 - 209622 x^2
         - 69166 x + 49036
\]
with $\#C(\Q) \ge 452$; the quotient is (at least)~$78.88$.


\section{Computations} \label{S:Data}

Our data come from several sources.

\subsection{Computations with (very) small curves}

This began as a project whose aim it was to decide, for every genus~2
curve $C \in \CC_3$, whether it possesses rational points. This experiment
is described in~\cite{BruinStollExp}, with more detailed explanation
of the various methods used in~\cite{BruinStoll2D,BruinStollMWS,BruinStollFG}.

These computations were later extended by the author. For those curves that
do have rational points, we proceeded to find all rational points, or at
least all rational points up to a height bound that is so large that we
can safely assume that no larger points exist.

More precisely, the following was done. We determined a generating set
for the Mordell-Weil group of the Jacobian of every curve (in a small
number of cases, the rank is not yet proved to be correct: there is a
difference of~$2$ between the rank of the known subgroup and the 2-Selmer
rank, which very likely comes from nontrivial elements of order~2 in
the Shafarevich-Tate group). When the Mordell-Weil rank~$r$ is zero,
the set $C(\Q)$ of rational points on~$C$ can be trivially determined.
When $r = 1$, a combination of Chabauty's method and the Mordell-Weil
sieve can be used to determine~$C(\Q)$; this is described in~\cite{BruinStollMWS}.
For $r = 2$, we can still use the Mordell-Weil sieve in order to find
all points up to a height of $H = 10^{1000}$ in reasonable time. For
$r > 2$, the sieving computation would take too long; in these cases,
we have used a lattice point enumeration procedure on the Mordell-Weil
group to find all points up to $H = 10^{100}$. The following table
summarizes what was done and gives the number of curves (up to isomorphism)
for each value of the rank~$r$. We denote the set of rational points
on~$C$ up to height~$H$ by $C(\Q)_H$.

\medskip

\begin{center}
\begin{tabular}{|c|r|l|} \hline
  $r = 0$ & 14\,010 curves & \strut $C(\Q)$ is determined. \\
  $r = 1$ & 46\,575 curves & \strut $C(\Q)$ is determined. \\
  $r = 2$ & 52\,227 curves & \strut $C(\Q)_H$ is determined for $H = 10^{1000}$. \\
  $r = 3$ & 22\,343 curves & \strut $C(\Q)_H$ is determined for $H = 10^{100}$. \\
  $r = 4$ &  2\,318 curves & \strut $C(\Q)_H$ is determined for $H = 10^{100}$. \\
  $r = 5$ &      17 curves & \strut $C(\Q)_H$ is determined for $H = 10^{100}$. \\\hline
\end{tabular}
\end{center}

\medskip

Under the reasonable assumption that there are no points on these curves
of height $> 10^{100}$ (note that the largest point that we found has
height about $2 \cdot 10^5$), plus assuming that all the ranks are correct,
this means that we have complete information on all rational points on
curves in~$\CC_3$. We plan to extend our computations to~$\CC_4$ eventually.

\subsection{All points with $H < 2^{14}$ on curves with $N \le 10$}

Since $N = 3$ is rather small, we also tried to get some information on
somewhat larger curves. The author has written a program {\tt ratpoints}
(see~\cite{Ratpoints} for a description) that uses a quadratic sieve and
fast bit-wise operations to search for rational points on hyperelliptic
curves. On current hardware, it takes about 10~ms on average to find all
points up to height $H = 2^{14}-1 = 16383$ on a genus~2 curve.

We used up to 20~machines from the CLAMV teaching lab at Jacobs
University Bremen for about one week in January~2008 to let {\tt ratpoints}
find all these points on all curves in~$\CC_{10}$. If $f \in \Z[x]$ is
the polynomial defining the curve, then it is only necessary to look at
one representative of the set
\[ \{f(x), f(-x), x^6 f(1/x), x^6 f(-1/x)\} \,, \]
since the corresponding curves are isomorphic and the isomorphism preserves
the height of the rational points. The total number of curves to be
considered was therefore roughly $21^7/4 \approx 450 \cdot 10^6$, for a
total of more than~$100$ CPU days (the average time per curve on these
machines was about 20~ms).

This gives us precise information on the frequency of points of
height $< 2^{14}$ on curves with $N \le 10$. It also gives us close
to complete information on curves with many points in this range, since
curves with many points seem to have reasonably small points. We might
have missed a few curves with (comparatively) many points that have one
(or more?) additional point pair(s).

We plan to extend these computations to $N \le 20$ (and possibly beyond),
with the same height bound, once we have suitable hardware at our disposal.

\subsection{Small curves with many points}

To get some more data on curves with many points, we conducted a systematic
search for curves in~$\CC_{50}$ with many points, making use of the
observation that all curves in~$\CC_{10}$ that have comparatively many points
tend to have points with $x$-coordinates $0$, $\infty$, $1$ and~$-1$.
Putting in the conditions that $F(0,1)$, $F(1,0)$, $F(1,1)$ and~$F(-1,1)$
have to be squares reduces the search space to a sufficient extent so
that a search up to $N = 50$ is possible. The point search was first done
with the bound $H = 2^{10} - 1 = 1023$; for those curves that had more than
a certain number of points in this range, points were then counted up to
height $H = 2^{17} - 1 = 131\,071$.

Based on the observation that all but one of the best curves that this
computation revealed also have rational points at $x = \pm 2$ (maybe after
a height-preserving isomorphism), we did a further systematic search for
curves in~$\CC_{200}$ having rational points at all $x \in \{\infty,0,1,-1,2,-2\}$.
Here the point search was done in three steps, using height bounds of
$2^{12} - 1 = 2047$, $2^{14} - 1 = 16383$, and finally $2^{17} - 1 = 131071$.
Two threshold values for the number of points were used in order to decide
whether to search for more points on a given curve.

We plan to extend these computations, too.

\subsection{Curves with many points in Elkies' families}

Noam Elkies was so kind to provide us with explicit formulas for five
ternary sextics $S(t,u,v)$ that admit many rational lines~$\ell$ on which $S$
restricts to a perfect square. Setting the restriction of~$S$ to a generic
line equal to a square gives a curve of genus~2 that has a pair of rational
points over each intersection point with a line~$\ell$ as above. In this way,
we obtain a 2-dimensional family of genus~2 curves with more than 50~pairs
of rational points. We have conducted a systematic search among all lines
$at + bu + cv = 0$ with $a,b,c \in \Z$ and $\max\{|a|,|b|,|c|\} \le 500$
in order to find curves with many points in these families.

There are two features of our computation that merit special mention.
The first is that we used as a preliminary selection step a product
$\prod_{p < X} \#C(\F_p)/p$ with $X = 200$, which was required to be
above a certain threshold value. The rationale behind this is that we
expect a curve with many rational points also to have more $\F_p$-points
than a random curve. Similar ideas have been used before. 
Note that each factor only depends on the reduction of the line $at+bu+cv=0$
mod~$p$, so that we can precompute the relevant values and reduce the
computation of the factors in the product to a table lookup.

The second is a systematic way of finding new rational points from known
ones. If there are five rational points $(x_i/z_i,y_i/z_i^3)$ on a genus~2 curve~$C$
that lie on a cubic $y = \alpha x^3 + \beta x^2 + \gamma x + \delta$, then
the sixth intersection point of this cubic with~$C$ is again a rational
point. The condition is equivalent to the vanishing of the determinant
of the matrix with rows $(z_i^3, x_i z_i^2, x_i^2 z_i, x_i^3, y_i)$.
For reasons of efficiency, we have written a C~program that first computes
these determinants mod~$2^{64}$ using native machine arithmetic; whenever
a determinant appears to be zero, this is checked using exact arithmetic,
and if the sixth intersection point is not yet known, it is recorded.
We have applied this procedure to the points of height up to~$10^5$ we found 
using {\tt ratpoints}. This can produce quite a number of additional points
of considerable height. For example, we were able to find eight more points
on Stahlke's curve from~\cite{Stahlke}, so that this curve must have at
least~374 rational points. Of course, in this way we can only find points within the
subgroup of the Mordell-Weil group generated by the known points.


\end{document}